\documentclass{amsart}
\usepackage{amssymb} 
\usepackage{amsmath} 
\usepackage{comment}
\usepackage[OT2,T1]{fontenc}
\DeclareSymbolFont{cyrletters}{OT2}{wncyr}{m}{n}
\DeclareMathSymbol{\Sha}{\mathalpha}{cyrletters}{"58}
\DeclareMathOperator{\spn}{Span}

\DeclareMathOperator{\ch}{char}
\DeclareMathOperator{\CH}{CH}
\DeclareMathOperator{\inv}{inv}
\DeclareMathOperator{\Br}{Br}

\DeclareMathOperator{\Pic}{Pic}

\newcommand{\up}{{\upsilon}}

\DeclareMathOperator{\ord}{ord}

\newcommand{\Sb}{S^{\mathrm{bd}}}
\newcommand{\Sg}{S^{\mathrm{gd}}}
\newcommand{\A}{\mathbb{A}}
\newcommand{\Q}{\mathbb{Q}}
\newcommand{\Z}{\mathbb{Z}}

\newcommand{\C}{\mathbb{C}}
\newcommand{\R}{\mathbb{R}}
\newcommand{\F}{\mathbb{F}}
\newcommand{\PP}{\mathbb{P}}

\newcommand{\cT}{\mathcal{T}}

\newcommand{\OO}{{\mathcal O}}

\newcommand{\uu}{\mathbf{u}}
\newcommand{\vv}{\mathbf{v}}

\newcommand{\xx}{\mathbf{x}}

\begin {document}

\newtheorem{thm}{Theorem}
\newtheorem*{thms}{Theorem}
\newtheorem{lem}{Lemma}[section]
\newtheorem{prop}[lem]{Proposition}

\newtheorem{cor}[lem]{Corollary}
\newtheorem{conj}[lem]{Conjecture}

\theoremstyle{definition}

\newtheorem*{ex}{Example}

\theoremstyle{remark}

\title[Mordell-Weil Generators for Cubic Surfaces]{
On the Number of 
Mordell-Weil Generators \\ for Cubic Surfaces
}
\author{Samir Siksek}
\address{Mathematics Institute\\
	University of Warwick\\
	Coventry\\
	CV4 7AL \\
	United Kingdom}

\email{s.siksek@warwick.ac.uk}
\date{\today}
\thanks{The author is supported by an EPSRC Leadership Fellowship}


\begin {abstract}
Let $S$ be a smooth cubic surface over a field $K$.
It is well-known that new $K$-rational points
may be obtained from old ones by secant and tangent constructions.
A Mordell-Weil generating set is a subset $B \subset S(K)$ of
minimal cardinality which generates $S(K)$
via successive secant and tangent constructions. Let $r(S,K)=\#B$.
Manin posed what is known as the Mordell-Weil problem for
cubic surfaces: if $K$ is finitely generated over its prime subfield
then $r(S,K)< \infty$. In this paper, we prove a special case of
this conjecture. Namely, if $S$ contains a pair of skew lines, both
defined over $K$, then $r(S,K)=1$. 

One of the difficulties in studying the secant and tangent process
on cubic surfaces is that it does not lead to an associative
binary operation as in the case of elliptic curves. As a partial
remedy we introduce an abelian group $H_S(K)$ 
associated to a cubic surface $S/K$,
naturally generated by the $K$-rational points, which retains much
information about the secant and tangent process. In particular,
$r(S,K)$ is large as soon as the minimal number of generators of $H_S(K)$
is large. In situations where $K$ is a global field, and weak
approximation holds for $K$-rational points on $S$, the group
$H_S^0(K)$ surjects onto $\prod_{\upsilon \in \Delta} H_S^0(K_\upsilon)$
for any finite set $\Delta$ of places, 
where $H_S^0$ is 
\lq the degree $0$ part of $H_S$\rq. 
We use this to construct 
a family of smooth cubic surfaces over $\Q$ 
such that $r(S,\Q)$ is unbounded in this family. 
This is the cubic surface analogue of the 
unboundedness of ranks conjecture for elliptic curves over $\Q$.
\end {abstract}
\maketitle

\section{Introduction}

Let $C$ be a smooth plane cubic curve over $\Q$. 
The Mordell-Weil Theorem
can be restated as follows: there is a finite subset $B$ of
$C(\Q)$ such that the whole of $C(\Q)$ can be obtained
from this subset by drawing secants and tangents through
pairs of previously constructed points and consecutively adding
their new intersection points with $C$.
It is conjectured that a minimal such $B$ can be arbitrarily large;
this is indeed the well-known conjecture that there are elliptic curves
with arbitrarily large ranks. This paper is concerned with
the cubic surface analogues of the Mordell-Weil Theorem and the unboundedness
of ranks.

Let $K$ be a field and let $S$ be a smooth cubic surface over $K$
in $\PP^3$. 
By a $K$-line we mean a line $\ell \subset \PP^3$ that is
defined over $K$. 
If $\ell \not\subset S$ then
$\ell \cdot S=P+Q+R$ where $P$, $Q$, $R\in S$. If
any two of $P$, $Q$, $R$ are $K$-points then so is the third.
The line $\ell$ is tangent at $P$ 
if and only if $P$ appears more than once in the sum $P+Q+R$.
If $B \subseteq S(K)$, we shall write $\spn(B)$
for the subset of $S(K)$ generated from $B$ by successive
secant and tangent constructions. 
More formally, we define a sequence
\[
B=B_0 \subseteq B_1 \subseteq B_2 \subseteq \cdots \subseteq S(K)
\]
as follows. We let $B_{n+1}$ be the set of points $R\in S(K)$
such that either $R \in B_n$,
or for some $K$-line $\ell \not \subset S$ we have
$\ell \cdot S=P+Q+R$ where $P$, $Q \in B_n$. Then
$\spn(B)= \cup B_n$. In view of the Mordell-Weil Theorem
for cubic curves it is natural to ask, for $K=\Q$
 say, if there is some finite
subset $B \subset S(K)$ such that $\spn(B)=S(K)$. As far as
we are aware, the possible existence of such an analogue of the Mordell-Weil
Theorem was first mentioned by Segre \cite[page 26]{Segre43} in 1943.
Manin
\cite[page 3]{Ma1} asks the same question for fields $K$ finitely
generated over their prime subfields.
He calls this \cite{Ma2} the Mordell-Weil problem for cubic
surfaces. 
The results of numerical experiments by Zagier 
(described by Manin in \cite{Ma2})
and Vioreanu \cite{V} lead different experts to different opinions
about the validity of this Mordell-Weil conjecture. 
However, Manin \cite{Ma2} and Kanevsky and Manin \cite{KM} prove
the existence of finite generating sets for rational points on
certain Zariski open subsets of rationally trivial
cubic surfaces with {\bf modified} composition operations induced by
birational maps $S \dashrightarrow \PP^2$.  
However, we are not aware of even a single example in the literature where
the existence of a finite set $B$ that generates $S(K)$ via
the secant and tangent process is proven.
In this paper we give a positive answer to a
special case of the Mordell-Weil problem.


\begin{thm}\label{thm:mw}
Let $K$ be a field with at least $13$ elements.
Let $S$ be a smooth cubic surface over $K$.
Suppose $S$ contains a pair of skew lines both defined over $K$.
Let $P\in S(K)$ be a point on either line that is not an
Eckardt point. 
 Then $\spn(P)=S(K)$.
\end{thm}
An {\em Eckardt point} is a point where three of the lines
contained in $S$ meet.
Note that Theorem~\ref{thm:mw} does not make the assumption that the
 field is finitely generated over its prime subfield! It is certainly
true that every smooth cubic surface over an algebraically closed
field is generated by a single element.
The theorem probably holds for most fields with fewer that $13$ elements,
but the proof is likely to involve a tedious examination
of special cases and we have not attempted it.
It is well-known that a cubic surface with a pair of skew $K$-lines
$\ell$, $\ell^\prime$
is birational to $\ell \times \ell^\prime$. The proof of Theorem~\ref{thm:mw}
(given in Section~\ref{sec:mw})
is an extension  of the proof that $S$
is birational to $\ell \times \ell^\prime$. 

Now let us write 
\[ 
r(S,K):=\min \{\#B : \text{$B \subseteq S(K)$ and $\spn(B)=S(K)$} \}.
\]
We are unable to show that $r(S,K)$ is finite for cubic
surfaces without a skew pair of $K$-lines. However,
in some cases we can bound $r(S,K)$ from below.
We use this to
show that $r(S,\Q)$
is arbitrarily large as $S$ varies among smooth cubic surfaces over $\Q$.
To do this we introduce
and study a simple analogue of the Picard group of an elliptic curve.
Let 
\[
G_S(K)= \bigoplus_{P \in S(K)} \Z\cdot P
\]
be the free abelian group generated by the $K$-rational points of $S$.
Let $G^\prime_S(K)$ be the subgroup generated by all three point sums 
$P+Q+R$ with $P$, $Q$, $R \in S(K)$ such that
\begin{enumerate}
\item[(i)] 
there is $K$-line $\ell$ not contained in $S$ with $\ell \cdot S=P+Q+R$, or
\item[(ii)] there is a $K$-line $\ell$ contained in $S$ such that
$P$, $Q$, $R \in \ell$.
\end{enumerate}
The {\em degree map}  $\deg : G_S(K) \rightarrow \Z$
is given by $\deg(\sum a_i P_i)=\sum a_i$.
Let 
\[
G^{\prime\prime}_S(K)=\{D \in G^\prime_S(K) : \deg(D)=0\}.
\]
Let $H_S(K):=G_S(K)/G^{\prime\prime}_S(K)$. If $P \in S(K)$ we denote the
image of $P$ in $H_S(K)$ by $[P]$. The degree map remains
well-defined on $H_S(K)$:
we let $\deg : H_S(K) \rightarrow \Z$ be given  by
$\deg(\sum a_i [P_i ])=\sum a_i$. 
We shall write 
\[
H^0_S(K)=\{ D \in  H_S(K) : \deg(D)=0\}.
\]
 If $S(K) \neq \emptyset$
then the degree homomorphism clearly induces
an isomorphism 
\[
H_S(K)/H^0_S(K) \cong \Z.
\]
The group $H_S(K)$ will allow us to study $r(S,K)$.

\begin{thm}\label{thm:two}
Let $p_1,\dots,p_s$ ($s\geq 1$) be distinct primes such that
\begin{enumerate}
\item[(a)] $p_i \equiv 1 \pmod{3}$,
\item[(b)] $2$ is a cube modulo $p_i$.
\end{enumerate}
Let $M=\prod p_i$ and 
let $S=S_M/\Q$ be the smooth cubic
surface given by
\begin{equation}\label{eqn:one}
S_M : x^3+y^3+z(z^2+M w^2)=0. 
\end{equation}
Write $H_S(\Q)[2]$ for the $2$-torsion subgroup of $H_S(\Q)$. Then
$H^0_S(\Q) = H_S(\Q)[2]$ and
\[
r(S,\Q) \geq \dim_{\F_2} H_S(\Q)[2] \geq 2s. 
\]
\end{thm}
A prime $p$ satisfies conditions (a) and (b) of the theorem
if and only if the polynomial $t^3-2$ has three roots modulo $p$.
By the Chebotar\"{e}v Density Theorem such primes 
form a set with Dirichlet density $1/6$ (see for example
\cite[pages 227--229]{Heilbronn}).
We thus see that $r(S_M,\Q)$ becomes arbitrarily large as $M$ varies. 
The cubic surface $S_M$ has precisely one $\Q$-rational
line, which is given by $x+y=z=0$.

Our second theorem concerns a family of diagonal cubic surfaces.
To obtain a similar result for diagonal cubic surfaces we need to
assume the following conjecture.
\begin{conj}\label{conj} (Colliot-Th\'{e}l\`{e}ne)
Let $X$ be smooth and proper geometrically rational surface
over a number field $K$. Then the Brauer-Manin obstruction
is the only one to weak approximation on $X$.
\end{conj}
This conjecture is stated as a question by Colliot-Th\'{e}l\`{e}ne and Sansuc in \cite[Section V]{angers}, but since then has been stated
as a conjecture by
Colliot-Th\'{e}l\`{e}ne \cite[page 319]{CT}.
Indeed for the proof of
Theorem~\ref{thm:two} we need the fact that $S_M$ satisfies weak approximation, but fortunately
this follows from a theorem of Salberger and Skorobogatov on degree $4$ del Pezzo surfaces.

\begin{thm}\label{thm:three}
Let $p_1,\dots,p_s$ ($s \geq 1$) be distinct primes $\equiv 1 \pmod{3}$.
Let $M=3\prod p_i$ and 
let $S=S^\prime_M/\Q$ be the smooth cubic
surface given by
\begin{equation}\label{eqn:two}
S^\prime_M : x^3+y^3+z^3+M w^3=0. 
\end{equation}
Assume that the Brauer-Manin obstruction is the only one to 
weak approximation for $S$.
Then
\[
r(S,\Q) \geq \dim_{\F_3} H^0_S(\Q)/3 H^0_S(\Q) \geq 2s. 
\]
\end{thm}
Theorem~\ref{thm:three} is less satisfactory than Theorem~\ref{thm:two}
in two obvious ways. The first is that it is conditional on the yet unproven (and probably very difficult)
conjecture concerning weak approximation. The second is that  
we do not know if $H_S^0(\Q)$ contains any elements of
infinite order or is merely a torsion group.

The paper is organized as follows. 
In Section~\ref{sec:class} we briefly review what we
need from the geometry of cubic surfaces. The proof 
of Theorem~\ref{thm:mw} is given in Section~\ref{sec:mw}.
In Section~\ref{sec:pre} we prove some useful results about $H_S(K)$ that
follow from the definition and the elementary geometry of cubic surfaces.
In particular, we show that 
$H^0_S(K) =  H_S(K)[2]$ if $S$ contains a $K$-line, and that
$H^0_S(K) =0$ if $S$ contains a pair of skew $K$-lines.
In Section~\ref{sec:local} we study $H_S(K)$ for local fields $K$.
In particular we show that $H^0_S(K)$ is finite and that the
map $S(K) \rightarrow H_S(K)$ given by $P \mapsto [P]$ is 
locally constant. In Section~\ref{sec:real} we show that
$H^0_S(\R)=0$ or $\Z/2\Z$ depending on whether $S(\R)$
has one or two connected components. In Section~\ref{sec:weak},
for $K$ a number field and $\Delta$ a finite set
of places of $K$, we study the diagonal map
$\mu_\Delta : H_S(K) \rightarrow \prod_{\upsilon \in \Delta} H_S(K_\upsilon)$;
for example if $S/K$ satisfies weak approximation then 
we show that the map is surjective. The proofs of Theorems~\ref{thm:two}
and~\ref{thm:three} essentially boil down to proving
(enough of) weak approximation and then estimating the size
of the target space of the diagonal map $\mu_\Delta$ 
for $\Delta=\{p_1,\dotsc,p_s\}$.
To this end we briefly introduce the Brauer-Manin obstruction
(Section~\ref{sec:BM}) and apply it in Sections~\ref{sec:BMSM} and~\ref{sec:BMSMd} to prove (enough of) weak approximation for the surfaces
$S_M$ and $S_M^\prime$, where in the latter case
we are forced to assume Colliot-Th\'{e}l\`{e}ne's Conjecture~\ref{conj}. 
Finally we must study $H_S(\Q_p)$ for the surfaces $S=S_M$,
$S=S_M^\prime$.
Let $C$ be the plane genus $1$ curve given by the equation
\begin{equation}\label{eqn:C}
C : x^3+y^3+z^3=0.
\end{equation}
For a prime $p \neq 3$ we shall denote $C_p=C \times \F_p$.
For $p =p_1,\dotsc,p_s$
it is easily seen that $S_M$ and $S_M^\prime$ both reduce to a cone over $C_p$.
This fact is crucial to the proofs of both Theorems~\ref{thm:two} and~\ref{thm:three}.
In Section~\ref{sec:C} we shall briefly study $\Pic^0(C_p)/2\Pic^0(C_p)$ and $\Pic^0(C_p)/3\Pic^0(C_p)$.
Section~\ref{sec:red} quickly reviews good choices of parametrizations of lines in $\PP^3$:
a good choice is one that still parametrizes the lines after reduction modulo
 $p$.
In Section~\ref{sec:SM}
 we shall construct, for the surface $S=S_M$,
 a surjective homomorphism $H^0_S(\Q_p) \rightarrow \Pic^0(C_p)/2\Pic^0(C_p)$.
In Section~\ref{sec:proofone} we use this and the surjectivity
of the diagonal map $\mu_\Delta$ for $\Delta=\{p_1,\dotsc,p_s\}$
to deduce Theorem~\ref{thm:two}.
We then turn our attention to the surface $S=S_M^\prime$. 
In Section~\ref{sec:SMd}
we construct a
surjective homomorphism $H^0_S(\Q_p) \rightarrow \Pic^0(C_p)/3\Pic^0(C_p)$, and conclude the proof
of Theorem~\ref{thm:three} in Section~\ref{sec:prooftwo}.

One reason why $H_S(K)$ may be of interest
is that it seems to be intimately related to
the Chow group $\CH_0(S)$ of zero-cycles on $S/K$. 
It is straightforward to see that elements of $G_S^{\prime\prime}(K)$
are zero-cycles that are rationally equivalent to $0$.
Thus we have a natural homomorphism
\[
\epsilon : H_S(K) \rightarrow \CH_0(S).
\]
In subsequent papers we plan to address the relationship of $H_S(K)$ with
$\CH_0(S)$, as well as various constructions found in \cite{Ma1}
such as Moufang loops, $R$-equivalence and universal equivalence. 
We will also make a more extensive study of
$H_S(K)$ for $K$ a finite or local field, 
and the natural pairing $\Br(S) \times \prod H_S(K_\upsilon) \rightarrow \Q/\Z$
induced by the corresponding pairing for the Chow group of zero cycles. 

I thank Felipe Voloch for drawing my attention to \cite{V}
(via the website {\tt mathoverflow.net}).
I am grateful to Martin Bright,  David Holmes, Miles Reid
and Damiano Testa
for helpful discussions, and the referee for many corrections
and comments.
I would like to thank Jean-Louis Colliot-Th\'{e}l\`{e}ne for 
useful correspondence
regarding the above conjecture and for drawing my attention to \cite{SS}.
In particular, Professor Colliot-Th\'{e}l\`{e}ne points out 
that it should be possible to deduce unboundedness results
for cubic surfaces 
similar to our Theorems~\ref{thm:two} and~\ref{thm:three}
from unboundedness results for the number of $R$-equivalence
classes of Ch\^{a}telet surfaces \cite[Theorem 8.13]{css}.


\section{Some Geometry} \label{sec:class}

We shall need some basic material on the geometry of cubic surfaces.
We do not claim any originality in this section, although we occasionally
sketch proofs of well-known statements in order to verify
that these hold in small positive characteristic. 

\begin{thm} (Cayley-Salmon) \label{thm:27}
Every non-singular cubic surface over an algebraically closed field contains exactly $27$ lines.

Every line $\ell$ on the surface meets exactly $10$ other lines,
which break up into $5$ pairs $\ell_i$, $\ell^\prime_i$ ($i=1,\dots,5$) such that
$\ell$, $\ell_i$ and $\ell_i^\prime$ are coplanar, 
and $(\ell_i \cup \ell_i^\prime) \cap (\ell_j \cup \ell_j^\prime)=\emptyset$ for $i \ne j$.
\end{thm}
\begin{proof}
For a proof see \cite[V.4]{Hartshorne} or \cite[Section IV.2]{Shaf}.
\end{proof}

For now $S$ will denote a smooth cubic surface in $\PP^3$ 
over a 
field $K$, defined by homogeneous cubic polynomial $F \in K[x_0,x_1,x_2,x_3]$.

For a point $P \in S(\overline{K})$,
 we shall denote the tangent plane to $S$ at $P$ by $\Pi_P$.
This is given by 
$\Pi_P : \nabla{F}(P) \cdot \xx=0$. 
We shall write $\Gamma_P$ for the plane curve $S \cdot \Pi_P$. It is
easy to check (using the smoothness of $S$)
that $\Gamma_P$ does not contain any multiple components.
It is a degree $3$ plane curve which is singular at $P$.
If $\Gamma_P$ is irreducible, it is nodal or cuspidal at $P$.
If $\Gamma_P$ is reducible then it is the union of a line and an irreducible conic,
or of three distinct lines.
\begin{lem}\label{lem:elem}
Let $P \in S(\overline{K})$. The curve $\Gamma_P$ contains every $\overline{K}$-line on $S$ that 
passes through $P$.
\end{lem}
\begin{proof}
By Euler's Homogeneous Function Theorem, $P\cdot \nabla{F}(P)=3 \cdot  F(P)=0$. 
The line $\ell$
has a parametrization of the form $s P+t \vv$ with $(s:t) \in \PP^1$. Thus the polynomial
$F(s P+t\vv)$
vanishes identically. However, coefficient of $t s^{2}$ in this
polynomial is $(\nabla{F})(P)\cdot \vv$. This shows that $\ell$
is also contained in $\Pi_P$, and hence in $\Gamma_p$.
\end{proof}

A $\overline{K}$-line $\ell$ is called an {\em asymptotic line} (c.f.\ \cite[Section 2]{Voloch})
at $P \in S(\overline{K})$ if $(\ell \cdot S)_P \geq 3$. 
As $S$ is a cubic surface, it is seen that for an asymptotic line $\ell$ at $P$,
either $(\ell \cdot S)_P=3$ or $\ell \subset S$.
The asymptotic lines at $P$ are contained in $\Pi_P$.
Any line contained in $S$ and passing through $P$ is an
asymptotic line through $P$. 
The number of distinct asymptotic $\overline{K}$-lines at $P$
is either $1$, $2$ or infinity.
If $S$ has either $1$ or 
infinitely many asymptotic lines at $P$ then we shall call $P$ a {\em parabolic} point. 
The case where there are infinitely many asymptotic lines  at $P$
is special: in this case $\Gamma_P$ 
decomposes as a union of three $\overline{K}$-lines lying on $S$
and so the point
$P$ is an {\em Eckardt} point. If $P$ is parabolic
but not Eckardt, the curve $\Gamma_P$ has a cusp
at $P$. If $P$ is non-parabolic, then $\Gamma_P$ has a node
at $P$. 

Next we suppose $P \in S(K)$. 
Then $\Pi_P$ and $\Gamma_P$ are defined over $K$.
Suppose $P$ is non-parabolic. 
Thus there are precisely two (distinct) asymptotic lines at $P$.
If these two lines are individually defined 
over $K$ then we say that $P$ is {\em $K$-hyperbolic}.
Otherwise the two lines are defined over a 
quadratic extension of $K$ and conjugate;
in this case we shall say that $P$ is {\em $K$-elliptic}. If $K=\R$
then the terms parabolic, hyperbolic and elliptic agree with their usual
meanings in differential geometry: they 
correspond to points where the Gaussian curvature
is respectively $=0$, $>0$ and $<0$.

We shall also need to study the number of parabolic points on a 
line lying on a cubic surface. Let ${\PP^3}^*$ be the dual projective
space and write $\gamma : S \rightarrow {\PP^3}^*$ for the 
{\em Gauss map} which sends a point to its tangent plane. 
A useful characterisation of parabolic points is that they
are the points of ramification of the Gauss map \cite[Section 2]{Voloch}.
If $\ell \subset S$ and $P \in \ell$, then $\ell$ is contained
in the tangent plane $\Pi_P$. The family of planes through $\ell$ can be
identified with $\PP^1$ and once such an identification is fixed we let
$\gamma_\ell : \ell \rightarrow \PP^1$ be the map that sends a point
on $\ell$ to its tangent plane through $\ell$.

\begin{lem}\label{lem:paraline}
Let $\ell$ be a $K$-line contained in $S$. Then
every $P\in \ell(K)$ is either parabolic or $K$-hyperbolic.
\begin{enumerate} 
\item[(i)] If $\ch(K) \neq 2$ then $\gamma_\ell$ is separable.
Precisely two points
$P \in \ell(\overline{K})$
are parabolic, and so there are at most two Eckardt points on $\ell$.
\item[(ii)] If $\ch(K)=2$ and $\gamma_\ell$ is separable then
there is precisely one point $P \in \ell(\overline{K})$
which is parabolic and so at most one Eckardt point on $\ell$.
\item[(iii)] If $\ch(K)=2$ and $\gamma_\ell$ is inseparable then
every point $P \in \ell(\overline{K})$ is parabolic
and the line $\ell$ contains exactly $5$ Eckardt points.
\end{enumerate}
\end{lem}
\begin{proof}
This is a well-known classical result for cubic surfaces over
the reals or complexes; see for example \cite[pages 103--104]{Segre}.
Let $P \in \ell(K)$. 
First we would like to show that
any $P \in \ell(K)$ is either parabolic or $K$-hyperbolic. Suppose $P$
is not parabolic. Then $S$ has precisely two asymptotic lines 
passing through $P$. 
One of these is $\ell$. Since this pair of asymptotic lines
must be $K$-rational as a whole, the other asymptotic line is $K$-rational.
Thus $P$ is $K$-hyperbolic.

Next we would like to count the number of parabolic points on $\ell$.
By a projective transformation defined over $K$ we may suppose that $\ell$ passes through the point $(0:0:0:1)$,
that the tangent plane at this point is $x_0=0$ 
and that the line is $x_0=x_1=0$. 
Then $F$ has the form
\[
F=x_0 Q + x_1 R
\]
with $Q$, $R$ are homogeneous quadratic forms in $K[x_0,x_1,x_2,x_3]$. 
Now $\gamma_\ell : \ell \rightarrow \PP^1$ can be written as
$P \mapsto (Q(P) : R(P))$. The fact that $S$ is non-singular implies
that $Q$ and $R$ do not simultaneously vanish along the line $\ell$.
Thus $\gamma_\ell$ has degree $2$. 
Suppose first that $\gamma_\ell$ is separable---this is always the case
if $\ch(K) \ne 2$. 
Applying 
the Hurwitz Theorem \cite[Section IV.2]{Hartshorne} to $\gamma_\ell$
immediately gives that the ramification divisor has degree $2$. 
If $\ch(K) \ne 2$ then the ramification is tame and so $\gamma_\ell$
is ramified at precisely two distinct $\overline{K}$-points. If $\ch(K)=2$,
then the ramification is wild and $\gamma_\ell$ is ramified at precisely one point.
Parts (i) and (ii) now follow as the parabolic points on $\ell$ are the ramification points
of $\gamma_\ell$, and as the Eckardt points are contained among the parabolic points.

Finally suppose that $\ch(K)=2$ and $\gamma_\ell$ is inseparable. Then
$\gamma_\ell$ is ramified at every point of $\ell$ and
so every point is parabolic. To complete the proof of (iii) we
must show that there are $5$ Eckardt points on $\ell$.
By Theorem~\ref{thm:27} there are $10$ lines on $S$ that
meet $\ell$. Let $\ell^\prime$ be one of these and let $P$ be
their point of intersection. Then $\ell$ and $\ell^\prime$
are distinct asymptotic lines to $S$ at $P$. The only way that $P$
can be parabolic is if there is a third line passing through $P$. 
Thus there are $5$ Eckardt points $\ell$ proving (iii).
%
\end{proof}

\begin{ex}
Many of the classical notions about cubic surfaces over $\C$
break down for cubic surfaces over a field of characteristic $2$. For example,
a cubic surface over $\C$ has $1$, $2$, $3$, $4$, $6$, $9$, $10$
or $18$ Eckardt points \cite[Section 100]{Segre} and any line
contains at most two Eckardt points. The following example shows that
this need not be the case in characteristic $2$.
Take $S$ to be the smooth cubic surface over $\F_2$ given by
\[
S:
x_0^2 x_2 + x_0^2 x_3 + x_0 x_1^2 + x_0 x_1 x_2 + x_0 x_3^2 + x_1^2 x_2 
+ x_1 x_2^2.
\]
The $27$ lines are rational over $\F_{64}$, and $S$
has $13$ Eckardt points. Three of the $27$ lines have $5$ Eckardt points
and the remaining all have exactly one.
\end{ex}

\section{Proof of Theorem~\ref{thm:mw}}\label{sec:mw}

Throughout this section we shall assume that $\#K\ge 13$.

\begin{lem}\label{lem:linegen}
Let $\ell$ be a $K$-line on $S$. Let $P \in \ell(K)$ be a point
that does not lie on any other line belonging to $S$. Then 
\[
\ell(K) \subseteq \Gamma_P(K) \subseteq \spn(P).
\]
\end{lem}
\begin{proof}
By Lemma~\ref{lem:elem}, 
we know $\ell \subseteq \Gamma_P$. Thus $\Gamma_P=\ell \cup C$
where $C$ is a conic contained in $\Pi_P$. Since $P$ does not
lie on any other line belonging to $S$, we know that $C$
is irreducible over $\overline{K}$. Now 
$P \in \ell \cap C$ since $P$ is either parabolic or $K$-hyperbolic by Lemma~\ref{lem:paraline}.
Thus
$\ell\cdot C = P+P^\prime$
where $P^\prime$ is also $K$-rational. 

Let $Q \in C(K)\backslash \{P,P^\prime\}$. Let $\ell^\prime$ be the
$K$-line joining $P$ and $Q$. Then $\ell^\prime \cdot S=2P+Q$
and so $Q \in \spn(P)$. Hence 
$C(K)\backslash \{P^\prime\} \subseteq \spn(P)$.

Now let $R \in \ell(K)\backslash \{P,P^\prime\}$. Let $Q_1 \in C(K)\backslash \{P,P^\prime\}$ and let $\ell^\prime$ be the $K$-line connecting $R$ and $Q_1$.
Then $\ell^\prime \cdot S=R+Q_1+Q_2$ where
$Q_2 \in C(K)\backslash \{P,P^\prime\}$. By the above, 
$Q_1$, $Q_2 \in \spn(P)$. Hence $R \in \spn(P)$. 
This shows that $\ell(K) \backslash \{P^\prime\}
\subseteq \spn(P)$. 

Finally we must show that $P^\prime \in \spn(P)$. 
There is nothing
to prove if $P =P^\prime$. Thus suppose $P \ne P^\prime$. Note that $\Pi_P=\Pi_{P^\prime}$.
By the proof of Lemma~\ref{lem:paraline},
the morphism $\gamma_\ell$ that sends points on $\ell$ to their tangent planes has degree $2$,
so for any $R \in \ell \backslash \{P,P^\prime\}$, $P^\prime$ 
will be a non-singular point of $\Gamma_R$.
The line $\ell$ meets precisely $10$ other lines lying on $S$ by Theorem~\ref{thm:27}.
The assumption that $\#K \ge 13$ forces the existence of $R \in \ell(K)$, different from $P$ and $P^\prime$,
and not lying on any line. Using the above argument with $R$ instead of $P$ we see that
$P^\prime \in \spn(R)$. However, $R \in \spn(P)$. This completes the proof.

\end{proof}

We now relax the hypotheses of Lemma~\ref{lem:linegen}
\begin{lem}\label{lem:linegen2}
Let $\ell$ be a $K$-line on $S$. Let $P \in \ell(K)$ 
and suppose that $P$ is not an Eckardt point. Then
\[
\ell(K) \subseteq \Gamma_P(K) \subseteq \spn(P).
\]
\end{lem}
\begin{proof}
To ease notation, we write $\ell_1$ for the line $\ell$.
If $P$ does not lie on any other line then this follows from Lemma~\ref{lem:linegen}.
Thus we may suppose that $P \in \ell_2$ where $\ell_1 \ne \ell_2$ but is not Eckardt. 
So $\ell_1$ and $\ell_2$ are distinct asymptotic lines at $P$, and so $P$ is non-parabolic.
By Lemma~\ref{lem:paraline},
the point $P$ must be $K$-hyperbolic, and so $\ell_2$ is defined over $K$. 
Now $\Gamma_P=\ell_1 \cup \ell_2 \cup \ell_3$ where $\ell_3 \subset S$
is a $K$-line not passing through $P$.
Write $P_{ij}=\ell_i \cdot \ell_j$; thus $P_{12}=P$ and the points $P_{ij}$ are distinct. 
Since the field is large enough, there is a point $Q \in \ell_3(K)$ such that $Q \ne P_{13}$, $P_{23}$,
and $Q$ does not lie on any other line contained in $S$. By Lemma~\ref{lem:linegen} we know that $\ell_3(K) \subseteq \spn(Q)$.
But if we let $\ell_Q$ be the line joining $Q$ with $P_{12}$, then $\ell_Q\cdot S=2P_{12}+Q$ and so $Q \in \spn(P_{12})$.
Thus $\ell_3(K) \subseteq \spn(P_{12})$. Similarly $\ell_1(K) \subseteq \spn(P_{23})$ and $\ell_2(K) \subseteq \spn(P_{13})$.
However, $P_{23}$, $P_{13} \in \ell_3(K)$. It follows that $\ell_i(K) \subseteq \spn(P_{12})$ for $i=1$, $2$, $3$.
This completes the proof. 
\end{proof}

\begin{lem}\label{lem:twolinegen}
Suppose
$S$ contains two skew
$K$-lines $\ell$, $\ell^\prime$.
Then 
\[
\ell^\prime(K) \subset \spn(\ell(K)).
\]
\end{lem}
\begin{proof}
As in Section~\ref{sec:class}, let $\gamma_\ell$ be the map that
sends a point $P \in \ell$ to the tangent plane $\Pi_P$ to 
$S$ at $P$. We know by the proof of Lemma~\ref{lem:paraline}
that $\gamma_\ell$ has degree $2$. Since $\ell \subset \Pi_P$
and $\ell$, $\ell^\prime$ are skew, we see that $\ell^\prime \cdot \Pi_P$
is a single point on $\ell^\prime$. Thus we can think of
$\gamma_\ell$ as a map $\ell \rightarrow \ell^\prime$
given by $P \mapsto \ell^\prime \cdot \Pi_P$.  

We claim the existence of some $P \in \ell(K)$ that is not
Eckardt such that $\gamma_\ell(P)$ is not Eckardt. Assume this
for the moment. It is clear that $\gamma_\ell(P) \in \Gamma_P(K)$.
By Lemma~\ref{lem:linegen2}, we see that $\gamma_\ell(P) \in \spn(P)$
and $\ell^\prime(K) \subseteq \spn(\gamma_\ell(P))$. 
Thus to complete the proof it is enough to establish our claim. 

By Lemma~\ref{lem:paraline}, $\ell^\prime$ contains at most $5$
Eckardt points and so it is enough to show that
\[
\# \gamma_\ell \left( \ell(K) \backslash \{ \text{Eckardt points}\} \right)
\geq 6.
\]
Suppose first that $\gamma_\ell$ is separable. The Eckardt points
on $\ell$ are at most $2$ by Lemma~\ref{lem:paraline}.
As $\gamma_\ell$ has degree $2$,
\[
\# \gamma_\ell \left( \ell(K) \backslash \{ \text{Eckardt points}\} \right)
\geq \frac{\#\ell(K) -2}{2} \geq 6.
\]
Finally suppose $\gamma_\ell$ is inseparable. Then 
$\gamma_\ell : \ell(K) \rightarrow \ell^\prime(K)$ is injective. 
Now there are $5$ Eckardt points on $\ell$ and
\[
\# \gamma_\ell \left( \ell(K) \backslash \{ \text{Eckardt points}\} \right)
\geq \#\ell(K) -5 \geq 9.
\]
\end{proof}

I am grateful to Damiano Testa for 
pointing out to me that a cubic surface containing skew $K$-lines $\ell$, $\ell^\prime$
is birational to $\ell\times \ell^\prime$ over $K$. This in essence is what the proof
of the following lemma is using.
\begin{lem}\label{lem:birattriv}
Suppose $S$ contains a pair of skew lines $\ell_1$ and $\ell_2$ both 
defined over $K$. Then
\[
\spn\left(\ell_1(K)\cup \ell_2(K) \right)=S(K).
\]
\end{lem}
\begin{proof}
Let $P$ be a $K$-point on $S$ not belonging to either line;
we will show that $P$ belongs to the span of $\ell_1(K) \cup \ell_2(K)$.
Let $\Pi_1$ be the unique plane containing $\ell_2$ and $P$,
and $\Pi_2$ the unique plane containing $\ell_1$ and $P$.
Since $\ell_1$ and $\ell_2$ are skew we know that $\ell_i \not \subset \Pi_i$.
Write $Q_i=\ell_i \cap \Pi_i$. 
Note that $P$, $Q_1$ and $Q_2$ are distinct points on $S$ that also belong
to the $K$-line $\ell=\Pi_1 \cap \Pi_2$. Suppose first that $\ell \not \subset S$. Then $\ell\cdot S=P+Q_1+Q_2$. Thus $P \in \spn\left(\ell_1(K)\cup \ell_2(K) \right)$ as required.

Next suppose that $\ell \subset S$. Then $\ell \subset \Gamma_{Q_1}$. 
If $Q_1$ is not Eckardt, then by
Lemma~\ref{lem:linegen2}, 
\[
P \in \ell(K)  \subseteq \Gamma_{Q_1}(K) \subseteq \spn(Q_1) \subseteq \spn(\ell_1(K)).
\]
Thus we may assume that $Q_1$ is Eckardt. Then $\Gamma_{Q_1}=\ell \cup \ell_1 \cup \ell_3$ where $\ell_3$ is also $K$-rational. Now, as $\ell_1$
and $\ell_2$ are skew, $\ell_2$ meets the plane $\Pi_{Q_1}$ in precisely
one point and this is $Q_2$. In particular, $\ell_2$ and $\ell_3$ are skew.
Thus $\ell_3(K) \subseteq \spn(\ell_2(K))$ by Lemma~\ref{lem:twolinegen}. Now let
$Q \in \ell_1(K) \backslash \{Q_1\}$. The line connecting $Q$ and $P$
meets $\ell_3$ in a $K$-point $R$, and so 
\[
P \in \spn(\ell_1(K) \cup \ell_3(K)) \subseteq \spn(\ell_1(K) \cup \ell_2(K)).
\]
\end{proof}

\begin{proof}[Proof of Theorem~\ref{thm:mw}]
This follows from Lemmas~\ref{lem:linegen2}, \ref{lem:twolinegen}
and \ref{lem:birattriv}.
\end{proof}

\section{Preliminaries on $H_S(K)$}\label{sec:pre}

We now begin our study of the group $H_S(K)$. Our
eventual aim is to prove Theorems~\ref{thm:two} and~\ref{thm:three},
but we will see other reasons why $H_S(K)$ is a
worthwhile object of study.
\begin{lem}\label{lem:line}
Let $C \subset S$ be either a $K$-line, a smooth plane
conic defined over $K$, or an irreducible plane cuspidal
or nodal cubic defined over $K$. Suppose $P$, $Q\in S(K)$
are points lying on $C$.
Moreover, in the cuspidal or nodal cubic case, suppose
that $P$, $Q$ are non-singular points on $C$.
Then $[P-Q]=0$ in $H_S(K)$.
\end{lem}

\begin{proof}
Suppose first that $C=\ell$ is a $K$-line lying on $S$
and $P$, $Q$ are $K$-points lying on $\ell$.
By definition of $G_S^\prime(K)$, we see that $3P$,  $2P+Q \in G_S^\prime(K)$.
Thus $[P-Q]=[3P-(2P+Q)]$ is zero in $H_S(K)$.

Next suppose that $C \subset S$ is a smooth plane conic
defined over $K$, and let $\Pi$
be the plane containing $C$. Then $S \cap \Pi=C \cup \ell$
where $\ell$ is a $K$-line. 
Let $P$, $Q$ be $K$-points on $C$.
The line joining $P$ and $Q$
meets $\ell$ in a $K$-point $P^\prime$. Likewise, the
tangent to $C$ at $Q$ meets $\ell$ in a $K$-point $Q^\prime$.
Thus $P+Q+P^\prime$, $2Q+Q^\prime \in G_S^\prime(K)$. Hence
\[
[P-Q]=[(P+Q+P^\prime)-(2Q+Q^\prime)]+[Q^\prime-P^\prime]=[Q^\prime-P^\prime].
\]
But by the previous part, $[Q^\prime]=[P^\prime]$ as both points
are on $\ell \subset S$. Hence $[P-Q]=0$.

Finally, suppose $C\subset S$ is a cuspidal or nodal plane cubic
defined over $K$, and let $R \in S(K)$ be the singular point on $C$.
Let $P$, $Q \in S(K)$ be non-singular points on $C$. By considering
the lines joining $P$ with $R$ and $Q$ with $R$, we see
that $P+2R$, $Q+2R \in G_S^\prime(K)$. Then $[P-Q]=0$.
\end{proof}

The reader will recall that whilst any rational point on
a cubic curve can be taken as the zero element, the description
of the group law is simpler and more pleasant if this point
is a flex. We shall now introduce a cubic surface analogue
of flexes.
Let $P \in S(K)$. We shall say that $P$ is {\em $K$-ternary} if
at least one of the asymptotic lines to $S$ through $P$ is
defined over $K$. 

\begin{lem}
If $P$ is $K$-hyperbolic then $P$ is $K$-ternary.
If $P$ is Eckardt then $P$ is $K$-ternary.
If $P$ is parabolic and $\ch(K) \ne 2$ then $P$ is $K$-ternary.
\end{lem}
\begin{proof}
The first part follows from the definition of $K$-hyperbolic.
 
Suppose $P$ is Eckardt. Of course each of the three lines through $P$
contained in $S$ is an asymptotic line. So we may suppose
that none of these is $K$-rational. Now let $\ell$ be any $K$-line
in $\Pi_P$ through $P$. Then $(\ell \cdot S)_P=3$ and $\ell$ is an asymptotic
$K$-line as required.

Finally, suppose $P$ is parabolic and $\ch(K) \ne 2$. We may suppose that $P$ is
non-Eckardt. Then $S$ has a unique asymptotic line $\ell$ through $P$. 
The asymptotic lines are obtained by solving a quadratic equation, and as there
is exactly one solution and $\ch(K) \ne 2$ the solution must be $K$-rational. 
\end{proof}

\begin{thm}\label{thm:terngen}
Suppose $P_0 \in S(K)$ is $K$-ternary.
Let $B\subseteq S(K)$ 
such that $\spn(B)=S(K)$.
Then the set $\{[P-P_0] : P \in B \}$ generates $H^0_S(K)$.
In particular, for any prime $p$, 
\[
r(S,K) \geq \dim_{\F_p} H^0_S(K)/p H^0_S(K).
\]
\end{thm}
\begin{proof}
Note that if $\ell$ is a $K$-rational line and $\ell\cdot S=P+Q+R$
with $P$, $Q$, $R \in S(K)$ then $P+Q+R$ and $3 P_0$ both belong
to $G_S^\prime (K)$ and so $[P-P_0]+[Q-P_0]+[R-P_0]=0$.
The first part of the lemma follows easily from this
and the fact that $H^0_S(K)$ is generated by elements of the form
$[P-P_0]$ with $P \in S(K)$.

For the second part, note that any set which generates $H^0_S(K)$
also generates the $\F_p$-vector space $H^0_S(K)/p H^0_S(K)$.
\end{proof}


\begin{thm} \label{thm:line}
Suppose $S$ contains a $K$-rational line $\ell$.
Then $H^0_S(K) = H_S(K)[2]$, where $H_S(K)[2]$ is the
$2$-torsion subgroup of $H_S(K)$.
\end{thm}
\begin{proof}
Let $Q$ be a $K$-rational point on $\ell$. 
We claim that $2P+Q \in G^\prime_K$
for all $P \in S(K)$. 
Assume our claim for a moment. Now $H^0_S(K)$
is generated by classes $[P^\prime - P]$ with 
$P$, $P^\prime \in S(K)$. 
By our claim, $2(P^\prime-P)=(2P^\prime+Q)-(2P+Q)$
is an element of degree $0$ in $G^\prime_S(K)$,
proving the theorem.

To prove our claim let $P \in S(K)$.
If $P \in \ell$ then $2P+Q \in G^\prime_S(K)$ by
definition of $G^\prime_S(K)$. Thus suppose $P \notin \ell$
and let $\Pi_P$ be the
tangent plane to $S$ at $P$. Now either $\Pi_P$ contains $\ell$
or it meets $\ell$ in precisely one $K$-rational point.
In either case, there is some $K$-rational point $Q^\prime \in \Pi_P \cap \ell$.
Since $P$ is not on $\ell$, we have $P \ne Q^\prime$.
Let $\ell^\prime$ be the unique $K$-line connecting $P$ and $Q^\prime$.
Now $\Gamma_P=\Pi_P \cap S$ is singular at $P$. 
Thus $\ell^\prime$ is tangent to $S$ at $P$,
and hence  $2P+Q^\prime \in G^\prime_S(K)$. 
However $Q-Q^\prime \in G^\prime_S(K)$ by Lemma~\ref{lem:line}.
This proves the claim.
\end{proof}

The proof of the following theorem is in essence a simplification of the
proof of Lemma~\ref{lem:birattriv}.
\begin{thm}\label{thm:twolines}
Suppose $S$ contains a pair of skew $K$-lines. Then $H_S^0(K)=0$.
\end{thm}
\begin{proof}
This can in fact be easily recovered from 
Theorem~\ref{thm:mw} if $\# K \geq 13$.
In any case, 
let $\ell_1$ and $\ell_2$ be a pair of skew $K$-lines contained in $S$.
Fix some $Q_i \in \ell_i(K)$.
We claim that $P+Q_1+Q_2 \in G_S^\prime(K)$ for every $P \in S(K)$. 
Let us assume this for the moment.
By the proof of Theorem~\ref{thm:line}, we know that
$2P+Q_2 \in G_S^\prime(K)$. Hence $P-Q_1 \in G_K^{\prime\prime}$;
in other words $P-Q_1 \equiv 0$ in $H_S^0(K)$. 
Since classes of the form $[P-Q_1]$ generate $H_S^0(K)$, we see that
$H_S^0(K)=0$.

It remains to prove our claim. Note, by Lemma~\ref{lem:line},
 any two points on
$\ell_i$ are equivalent. Moreover, if $P \in \ell_1(K)$ then
all we have to show is that $2Q_1+Q_2 \in G_S^\prime(K)$ which
is true from the proof of Theorem~\ref{thm:line}.
Thus we may suppose that $P \notin \ell_i$, $i=1$, $2$.
Let $\Pi_1$ be the unique plane containing $\ell_2$ and $P$,
and $\Pi_2$ the unique plane containing $\ell_1$ and $P$.
Since $\ell_1$ and $\ell_2$ are skew we know that $\ell_i \not \subset \Pi_i$.
Write $Q_i=\ell_i \cap \Pi_i$. 
Note that $P$, $Q_1$ and $Q_2$ are distinct points on $S$ that also belong to
to the $K$-line $\Pi_1 \cap \Pi_2$. By the definition of $G_S^\prime(K)$ we have
$P+Q_1+Q_2 \in G^\prime_S(K)$.
\end{proof}

\begin{cor} \label{cor:ac}
Let $K$ be algebraically closed. Then
$H^0_S(K)=0$.
\end{cor}

\section{$H_S(K)$ for Local Fields $K$}\label{sec:local}
In this section we let $K$ be a local field. This means
that $K=\R$ or $\C$, or that $K$ is a finite extension of $\Q_p$
for some prime $p$.
The purpose of this section is to prove the following theorem.
\begin{thm}\label{thm:local}
Let $K$ be a local field 
and let $S$ be a smooth cubic surface over $K$. Then
the following hold.
\begin{enumerate}
\item[(i)] The group
 $H_S^0(K)$ is finite.
\item[(ii)] The natural map $S(K) \rightarrow H_S(K)$ given by $P \mapsto [P]$
is locally constant, where $S(K)$ has the topology induced by $K$.
\end{enumerate}
\end{thm}

Of course, if $S(K)$ is empty there is nothing to prove,
so we shall suppose $S(K) \neq \emptyset$. 
Then $S(K)$ is a $2$-dimensional
$K$-manifold, and in particular has $K$-points not lying on the $27$ lines.
Such a $K$-point is called a {\em general point}.

The proof of Theorem~\ref{thm:local} 
requires a series of lemmas. The key lemma is the
following.

\begin{lem} \label{lem:key}
There exists some non-empty $W \subset S(K)$, open in
the topology induced by $K$, such that for all 
$P_1$, $P_2 \in W$ we have $[P_1]=[P_2]$
in $H_S(K)$.
\end{lem}
\begin{proof}
It is well-known
that the existence of a general $K$-point makes $S$ unirational,
and makes the set $S(K)$ dense in the Zariski topology. 
One proof of this is found in Manin's book~\cite[Chapter II, Theorem 12.11]{Ma1}.
Our lemma follows from a modification, given below, of Manin's proof.

Let $R$ be a general $K$-point on $S$.
Observe that the plane cubic curve $\Gamma_R$ is irreducible,
defined over $K$
and singular only at $R$.
Let $\cT_S$ be the projectivized tangent bundle of $S$. 
In other words, $\cT_S$ paramet\-rizes pairs $(Q,\ell)$ where
$Q$ is a point on $S$ and $\ell$ is a line tangent to $S$
at $Q$. As the bundle $\cT_S$ is locally trivial, 
there is a Zariski open subset $U \subset S$ containing $R$
and a local isomorphism $U \times \PP^1 \rightarrow \cT_S$ 
defined over $K$. Remove from $U$ the $27$ lines on $S$. 
Then the local isomorphism induces a morphism $\phi : U \times \PP^1 \rightarrow S$
defined over $K$ as follows: if $Q\in U$, $\alpha \in \PP^1$
and $\ell_\alpha$
is the line tangent to $S$ at $Q$ corresponding to $\alpha$
 then let $\phi(Q,\alpha)$ be the third point of intersection of
$\ell$ with $S$. 
Consider the restriction
 of $\phi$ to $(\Gamma_R \cap U) \times \PP^1$. The image
must be irreducible and contains $\Gamma_Q$ for all $Q \in \Gamma_R \cap U$
(including $\Gamma_R$).
If the image is $1$-dimensional then
$\Gamma_Q =\Gamma_R$ which is impossible for $Q \ne R$ since
$Q$ is non-singular in $\Gamma_R$ but singular in $\Gamma_Q$.
This shows that $(\Gamma_R \cap U) \times \PP^1 \rightarrow S$
is dominant. 

Let $V$ be a non-empty subset of   $\Gamma_R(K) \cap U(K)$, open
in the topology induced by $K$, such that $R \notin V$. 
Then $\phi(V \times \PP^1(K))$ contains a non-empty 
subset $W$ open in $S(K)$. To complete the proof it is enough to show that
the map $P \mapsto [P]$ is constant on $\phi(V \times \PP^1(K))$.
Suppose $P_1$, $P_2$ are in this set. So $P_i+2 Q_i = \ell_i \cdot S$
for some $K$-lines $\ell_1$, $\ell_2$ and for some $Q_1$, $Q_2 \in V$.
Now $Q_1$, $Q_2 \in \Gamma_R$, but neither is equal to $R$. Thus $[Q_1]=[Q_2]$
by Lemma~\ref{lem:line}. It follows that $[P_1]=[P_2]$ as desired.
\end{proof}

\begin{lem}\label{lem:local}
The map $S(K) \rightarrow H_S(K)$ given by $P \mapsto [P]$
is locally constant.
\end{lem}
\begin{proof}
By Lemma~\ref{lem:key}, there is some non-empty open $W \subset S(K)$
on which the map $P \mapsto [P]$ is constant.
We shall prove the lemma by using secant operations
to cover $S(K)$ by \lq translates\rq\ of $W$. 

Let $P_0 \in S(K) \backslash W$. Choose $Q_0 \in W$
such that the $K$-line $\ell$ joining $Q_0$ and $P_0$ is not
contained in $S$ and is not tangent to $S$ at either point.
Write $R$ for the $K$-point such that $\ell \cdot S=R+Q_0+P_0$.
Consider the rational map
\[
t_R : S \dashrightarrow S, \qquad
\text{$P \mapsto Q$ if $P$, $Q$ and $R$ are collinear}. 
\]
Clearly $t_R$ restricts to a local homeomorphism
of neighbourhoods of $P_0$ and $Q_0$. Thus
there is some open $V$ containing $P_0$, contained
in the domain of this local homeomorphism, such that $t_R(V) \subset W$.
Now let $P \in V$. 
Then $Q_0+P_0+R \in G_S^\prime(K)$ and $t_R(P)+P+R \in G_S^\prime(K)$.
Thus $[Q_0-t_R(P)]+[P_0-P]=0$. But $[t_R(P)]=[Q_0]$ as both
$t_R(P)$ and $Q_0$ are in $W$. Thus $[P]=[P_0]$. Hence $P \mapsto [P]$
is constant on the open neighbourhood $V$ of $P_0$.
\end{proof}

\begin{proof}[Proof of Theorem~\ref{thm:local}]
By Lemma~\ref{lem:local} there is a covering of 
$S(K)$ by open sets $U$ such that the function $P \mapsto [P]$ 
is constant when restricted to $U$. As $S(K)$ is compact,
we can assume that we have finitely many such $U$, say
$U_1,\dotsc,U_m$. Fix $P_i \in U_i$ for $i=1,\dotsc,m$.
The abelian group $H_S^0(K)$ is then
generated by the differences $[P_i-P_j]$.

It will be sufficient to show that for any two
points $P$, $Q \in S(K)$ 
the class $[P-Q]$ has finite order.  Let $\ell$ be the $K$-line
joining $P$ and $Q$. If $\ell \subset S$ then $[P-Q]=0$
by Lemma~\ref{lem:line}, so suppose that $\ell \not \subset S$.
Let $\Pi$ be a $K$-plane through $\ell$ that misses all the lines
on $S$ and such that 
the irreducible plane cubic $E =S \cap \Pi$ is non-singular at $P$ and $Q$.

If $E$ is singular, 
then $[P-Q]=0$ by Lemma~\ref{lem:line}. 
Thus we may suppose that $E$ is non-singular. Consider the
elliptic curve $(E,Q)$. For an integer $n$, there is 
a unique $R_n \in E(K)$ such that $n(P-Q) \sim R_n-Q$, where $\sim$
denotes linear equivalence on $E$. Now the group operations
on $(E,Q)$ are given by secants and tangents, so we know that
$n[P-Q]=[R_n-Q]$. However, from the properties of elliptic curves
over local fields, we may choose values $n \neq 0$ so that $R_n$ is arbitrarily
close to $Q$. But we already know that if $R_n$ is sufficiently
close to $Q$ then $[R_n]=[Q]$. 
This completes the proof.
\end{proof}

\section{$H_S(\R)$}\label{sec:real}
Let $S$ be a smooth cubic surface defined over $\R$.
We quickly summarise some well-known facts; for 
details see \cite{PT} and \cite{SJS}.
The cubic surface $S$ contains
3, 7, 15 or 27 real lines, and 
$S(\R)$ has either one or two connected components.
If $S(\R)$ has one connected component then this
component is non-convex. If $S(\R)$ has two connected components
then one is convex and the other non-convex. The non-convex
component contains all the real lines.

\begin{thm}\label{thm:real} 
If $S(\R)$ consists of one connected component then
$H_S^0(\R)=0$. If $S(\R)$ consists of two connected components then
$H_S^0(\R) \cong \Z/2\Z$. This isomorphism may be given 
explicitly as follows. Let $P_0$ be a point on the 
non-convex component. The map 
\[
[P-P_0] \mapsto \begin{cases}
\overline{0} & \text{if $P$ is in the non-convex component} \\
\overline{1} & \text{if $P$ is in the convex component}\\
\end{cases}
\]
extends to an isomorphism $H_S^0(\R) \cong \Z/2\Z$.
\end{thm}
\begin{proof}
By Theorem~\ref{thm:local}, the map $P \mapsto [P]$
is locally constant, and hence constant on each connected
component. Now $H_S^0(\R)$ is generated by
classes of differences $[P-Q]$, and so if there is only one component
then $H_S^0(\R)=0$.

Now suppose $S(\R)$ has two components. 
To see that the map given in the theorem extends to
an isomorphism $H_S^0(\R) \cong \Z/2\Z$,
it is enough to observe that if $P$, $Q$, $R$
are collinear real points on $S$, then either
all three are on the non-convex component,
or precisely one is on the non-convex component.
\end{proof}

\bigskip

\noindent{\bf Remark.} I am grateful to the referee
for the following remark. By a result of
Colliot-Th\'el\`ene and Ischebeck \cite{CI},
the degree $0$ part of $\CH_0(S/\R)$ is
isomorphic $(\Z/2\Z)^{s-1}$, where
$s$ is the number of real components. It follows
from this and Theorem~\ref{thm:real}
that the natural map $H_S(\R) \rightarrow \CH_0(S/\R)$
is an isomorphism.

\section{Weak Approximation and $H_S$}\label{sec:weak}

In this section $K$ denotes a number field and $\Omega$
the set of places of $K$. 
Denote the ad\`{e}les of $K$ by $\A_K$.
As usual, $S$ is a smooth cubic surface
over $K$, but we shall further suppose that $S(\A_K) \ne \emptyset$.
We say $S$ satisfies {\em weak approximation}
if the image of $S(K)$ in $S(\A_K)$ is dense. More generally,
let $\Sigma$ be a finite subset of $\Omega$,
and denote by $\A_K^\Sigma$ the ad\`{e}les of
$K$ with the $\Sigma$-components removed.
We say that $S$ satisfies {\em weak approximation}
away from $\Sigma$ if the image of $S(K)$ in $S(\A_K^\Sigma)$
is dense. If we assume Colliot-Th\'{e}l\`{e}ne's Conjecture~\ref{conj},
then it is easy in any given case to write down a finite set of places
 $\Sigma$ such that $S$ satisfies weak approximation away from $\Sigma$.
We shall not do this in general, but only for the surfaces
$S_M$ and $S_M^\prime$ in Theorems~\ref{thm:two} and~\ref{thm:three}.

\begin{thm}\label{thm:surject}
Let $K$ be a number field and $\Omega$ its places. 
Let $\Sigma$ be
a subset of $\Omega$ and let $\A_K^\Sigma$ denote the
ad\`{e}les of $K$ with the $\Sigma$ components removed.
Suppose that the image of $S(K)$ in $S(\A_K^\Sigma)$
is dense. Let $\Delta$ be a finite subset of $\Omega \backslash \Sigma$.
Then the diagonal map
\begin{equation}\label{eqn:diag}
H^0_S(K) \rightarrow \prod_{\upsilon \in \Delta} H^0_S(K_\upsilon)
\end{equation}
is surjective.
\end{thm}
\begin{proof}
The target space of the homomorphism in \eqref{eqn:diag} is generated by elements of the form $\left([P_\upsilon-Q_\upsilon]\right)_{\upsilon\in \Delta}$, thus it is enough to show that
such elements are in the image.
By the hypotheses, $S(K)$ is dense in $\prod_{\upsilon \in \Delta} S(K_\upsilon)$. Recall, from Theorem~\ref{thm:local}, that the maps $R \rightarrow [R]$ are
locally constant on the $S(K_\upsilon)$. Thus choosing $P$, $Q \in S(K)$
that sufficiently 
approximate $(P_\upsilon)_{\upsilon\in \Delta}$ and 
$(Q_\upsilon)_{\upsilon\in \Delta}$ will give an element
$[P-Q] \in H_S^0(K)$ whose image under \eqref{eqn:diag}
is $\left([P_\upsilon-Q_\upsilon]\right)_{\upsilon\in \Delta}$.
\end{proof}

At first sight it seems that this theorem enables
us to disprove the Mordell-Weil conjecture for cubic
surfaces
simply by taking the set $\Delta$ to be arbitrarily large
and forcing $H_S^0(K)$ to surject onto larger and larger
groups. However, 
extensive---though not systematic---experimentation with cubic
surfaces, which we do not describe here, suggests that
$H^0_S(K_\upsilon)=0$ if $\upsilon$ is a non-Archimedean
place of good reduction for $S$ and $\upsilon \nmid 2$. 
By a place of good reduction we mean a non-Archimedean
place $\upsilon$ such that the polynomial defining $S$
has $\upsilon$-integral coefficients, and the reduction
of that polynomial modulo $\upsilon$ defines a smooth
cubic surface over the residue field.
We remark that for any place $\upsilon$ of good reduction
for $S$, it is known (e.g.\ \cite[Theorem A]{CT2}) that 
the degree $0$ part of the Chow group vanishes for $S \times K_\upsilon$.

We shall use Theorem~\ref{thm:surject} to
prove Theorems~\ref{thm:two} and~\ref{thm:three}. 
The first step is to prove weak approximation for
the surfaces $S_M/\Q$ in Theorem~\ref{thm:two}
and weak approximation away from $\{3\}$ for
the surfaces $S_M^\prime/\Q$ in Theorem~\ref{thm:three};
in the latter case our result will be conditional
on Colliot-Th\'{e}l\`{e}ne's Conjecture~\ref{conj}.
To prove weak approximation we introduce the Brauer-Manin obstruction
and study it for the surfaces $S_M$ and $S_M^\prime$.

\section{Brauer-Manin Obstruction: A Brief Overview}\label{sec:BM}

To proceed further we need to recall some
facts about the Brauer-Manin obstruction; 
for fuller details see \cite[Section 5.2]{Sk}. 
We continue with the notation of the previous section:
 $K$ is a number field, $\A_K$ the ad\`{e}les of $K$
and $\Omega$ the set of places of $K$.
The {\em Hasse reciprocity law} states that the following
sequence of abelian groups is exact:
\[
0 \rightarrow \Br(K) \rightarrow \sum_{\up \in \Omega} \Br(K_\up) \rightarrow \Q/\Z \rightarrow 0.
\]
Here the third map is the sum of local invariants 
$\inv_\up : \Br(K_\up) \hookrightarrow \Q/\Z$. 
Let $X$ be a smooth, projective and geometrically integral variety 
over a number field $K$.
Let $\Br(X)$ be the Brauer group of $X$
and denote by $\Br_0(X)$ the image of $\Br(K)$ in $\Br(X)$. 
Consider the pairing
\[
\langle~,~\rangle : \Br(X) \times  X(\A_K) \rightarrow \Q/\Z,
\qquad 
\langle A, ( P_\up) \rangle =\sum_{\up \in \Omega} \inv_\up (A(P_\up)).
\]
This is the {\em adelic Brauer-Manin pairing} and 
satisfies the following properties.
\begin{enumerate}
\item[(i)] If $A \in \Br_0(X) \subset \Br(X)$ and $(P_\up) \in X(\A_K)$
then $\langle A , (P_\up) \rangle=0$.
\item[(ii)] If $P \in X(K)$ then $\langle A , P \rangle=0$
for every $A \in \Br(X)$.
\item[(ii)] For any $A \in \Br(X)$, the map 
\[
 X(\A_K) \rightarrow \Q/\Z, \qquad
(P_\up) \mapsto \langle A , (P_\up) \rangle 
\]
 is continuous where  
 $\Q/\Z$ is given the discrete topology.
\end{enumerate}
We define
\[
X(\A_K)^{\Br(X)} = 
\{ (P_\up) \in  X(\A_K) : 
\text{ $\langle A , (P_\up) \rangle=0$ for all $A \in \Br(X)/\Br_0(X)$}  \}.
\]
By the above we know that
\[
\overline{X(K)} \subseteq X(\A_K)^{\Br(X)},
\]
where $\overline{X(K)}$ is the closure of $X(K)$.
We say that
{\em the Brauer-Manin obstruction is the only obstruction to
weak approximation} if $\overline{X(K)}=X(\A_k)^{\Br(X)}$. 

\section{The Brauer-Manin Obstruction for $S_M$}\label{sec:BMSM}
In this section we prove the following proposition.
\begin{prop}\label{prop:dense}
Let $p_1,\dots,p_s$ ($s \geq 1$) and $M$ be as in Theorem~\ref{thm:two}.
Let $S=S_M/\Q$ be the cubic surface given by $\eqref{eqn:one}$.
Then $S$ satisfies weak approximation. In particular,
the homomorphism 
\[
H^0_{S}(\Q) \rightarrow \prod_{p=p_1}^{p=p_s}  H^0_S(\Q_p)
\]
is surjective.
\end{prop}
To prove this proposition we shall need a theorem of
Salberger and Skorobogatov on del Pezzo surfaces of degree $4$.
The cubic surface $S_M$ is birational to a degree $4$ del Pezzo $X$
given by
the following smooth intersection of two quadrics in $\PP^4$:
\[
X_M : \begin{cases}
x^2-xy+y^2+zt=0, \\
z^2+Mw^2 -xt-yt=0.
\end{cases}
\] 
The map $X_M \dashrightarrow S_M$ is the obvious one $(x,y,z,w,t) \mapsto (x,y,z,w)$.
\begin{lem}
$\Br(X_M)/\Br_0(X_M)$ is trivial.
\end{lem}
\begin{proof}
To determine $\Br(X_M)/\Br_0(X_M)$ we can
use the recipe in \cite{Swd1}
or the more detailed recipe in \cite{BBFL}. The first step is 
to write down the $16$ lines on $X_M$. We did this by writing down
and solving the equations for the corresponding zero-dimensional 
Fano scheme (see for example \cite[Section IV.3]{EH}).
Write $\theta=\sqrt[3]{2}$ and let $\zeta$ be a primitive cube root of unity.
There are two Galois orbits of lines. The first orbit has
four lines
and a representative is
\[
L_1 :
\begin{cases}
x+\zeta y=0, \\
z+\sqrt{-M} w=0, \\
t=0.\\
\end{cases}
\]
The second orbit has 12 lines and a representative is
\[
L_2 :
\begin{cases}
-\theta z+ \theta \sqrt{-M}w + t=0,\\ 
3 \theta^2 x + (2\zeta - 2)\theta z + (\zeta + 2) t=0, \\
3 \theta^2 y + (-2 \zeta - 4) \theta z + (-\zeta + 1)  t=0.\\
\end{cases}
\]
From the size of these orbits we know \cite[Proposition 13]{BBFL} 
 that $\Br(X_M)/\Br_0(X_M)$ is trivial. 
\end{proof}

We shall also need the following theorem.
\begin{thms} (Salberger and Skorobogatov \cite[Theorem 6.5]{SS})
Let $K$ be a number field and $X$ a del Pezzo surface of 
degree $4$ over $K$ containing a rational point. 
Then $X(K)$ is dense in $ X(\A_K)^{\Br(X)}$.
\end{thms}

\begin{proof}[Proof of Proposition~\ref{prop:dense}]
By the above theorem of Salberger and Skorobogatov we know that
$X_M$ satisfies weak approximation. Now $S_M$ is birational to $X_M$,
and so by \cite[Lemma 5.5(c)]{SS} also satisfies weak approximation.
The last part of the proposition follows from Theorem~\ref{thm:surject}.
\end{proof}

\section{The Brauer-Manin Obstruction for $S_M^\prime$}\label{sec:BMSMd}

In this section we prove the following proposition.
\begin{prop}\label{prop:densed}
Let $p_1,\dots,p_s$ and $M$ be as in Theorem~\ref{thm:three}.
Let $\Sigma=\{3\}$. 
Let $S=S_M^\prime/\Q$ be the cubic surface in \eqref{eqn:two} and suppose
that the Brauer-Manin obstruction is the only obstruction
to weak approximation on $S$. Then $S$
satisfies weak approximation away from $\{3\}$.
In particular, the homomorphism
\[
H_S^0(\Q) \rightarrow \prod_{p=p_1}^{p=p_s} H_S^0(\Q_p)
\]
is surjective.
\end{prop}
\begin{proof} 
All the results we need for this proof are due
to Colliot-Th\'{e}l\`{e}ne, Kanevsky and Sansuc \cite{CKS},
though Jahnel's Habilitation  
summarizes these results in one convenient theorem \cite[Chapter III, Theorem 6.4]{Jahnel}. 
Indeed, we know that 
\begin{enumerate}
\item[(i)] $\Br(S)/\Br_0(S) \cong \Z/3\Z$. 
Fix $A \in \Br(S)$ that represents a non-trivial
coset of $\Br(S)/\Br_0(S)$. 
\item[(ii)] The image of 
\[
\langle~,~\rangle : \Br(S) \times S(\A_\Q) \rightarrow \Q/\Z
\]
is $\frac{1}{3} \Z/\Z$.
\item[(iii)] The map 
\[
S(\Q_p) \rightarrow \frac{1}{3} \Z/\Z, \qquad P \mapsto \inv_p(A,P)
\]
is surjective for all $p \mid M$.
\end{enumerate}
Now the strategy is clear. Suppose that 
$\mathcal{P}=(P_\up) \in S(\A_\Q^\Sigma)$.
Choose $P_3 \in S(\Q_3)$ such that
\[
\inv_3(A,P_3)=-\sum_{\up \neq 3} \inv_\up (A,P_\up).
\]
Let $\mathcal{P}^\prime \in S(\A_\Q)$ be the point obtained
from $\mathcal{P}$ by taking $P_3$ to be the component at $3$.
Then $\langle A, \mathcal{P}^\prime \rangle=0$. Since $A$
generates $\Br(S)/\Br_0(S)$ we know that 
$\mathcal{P}^\prime \in S(\A_\Q)^{\Br(S)}$. By our assumption that
the Brauer-Manin obstruction is the only one to weak approximation
we have that $\mathcal{P}^\prime$ is in the closure of the rational
points in $S(\A_\Q)$. The proposition follows.
\end{proof}

\section{A Plane Cubic}\label{sec:C}
We shall need to study the reduction of the cubic surfaces
\eqref{eqn:one} and \eqref{eqn:two} at the primes $p \mid M$.
We note that modulo $p$ they both reduce to cones over the plane
cubic curve $x^3+y^3+z^3=0$. In this section we collect some
information we need regarding the Picard group of this
cubic curve.

In this section $K$ is a field of characteristic $\neq 3$. 
Throughout the rest of the paper, $C/\Q$ will denote the
plane genus $1$ curve given by \eqref{eqn:C}.
Note that with the restriction imposed on the characteristic, $C \times K$ is smooth.
Let $\OO=(1:-1:0) \in C$. 
We note that $\OO$ is a flex
and so $C$ can be put into Weierstrass form by a projective
transformation that sends $\OO$ to the point at $\infty$
(see for example \cite[Proposition 2.14]{Knapp}).
It follows for $P$, $Q$, $R \in C(K)$ that
$P+Q+R-3\OO=0$ in $\Pic^0(C\times K)$ if and only if there is $K$-line $\ell$
in $\PP^2$ such that $\ell \cdot C=P+Q+R$.

\begin{lem}\label{lem:mod3}
Let $p \equiv 1 \pmod{3}$ be a prime and write $C_p$ for $C \times \F_p$.
Then 
\[
\dim_{\F_3} \Pic^0(C_p)/3\Pic^0(C_p)=2. 
\]
Moreover, 
each of the $9$ elements of $\Pic^0(C_p)/3\Pic^0(C_p)$ 
can 
be represented by the class of $P-\OO$ for some $P \in C(\F_p)$.
\end{lem} 
\begin{proof}
$\Pic^0(C_p)$ is a finite abelian group
isomorphic to $E(\F_p)$ where $E$ is the elliptic curve $(C,\OO)$.
Hence we have isomorphisms
\[
\Pic^0(C_p)/3\Pic^0(C_p) \cong \Pic^0(C_p)[3] \cong E(\F_p)[3],
\]
the first of which is of course non-canonical.
The assumption that $p \equiv 1 \pmod{3}$ ensures that the
nine flex points of $C$ are defined over $\F_p$; these
have the form $(1 : -\omega :0)$ and permutations of these
coordinates, where $\omega^3=1$. Now the $3$-torsion in
$\Pic^0(C_p)$ consists of the nine classes $[\OO^\prime-\OO]$
where $\OO^\prime$ is a flex point. Hence
\[
\Pic^0(C_p)/3\Pic^0(C_p) \cong \Z/3\Z \oplus \Z/3\Z.
\]
This proves the first part of the lemma.
The second part of the lemma follows since every element of $\Pic^0(C_p)$
is the class of $P-\OO$ for some $P\in C(\F_p)$.
\end{proof}

\begin{lem}\label{lem:mod2}
Let $p \equiv 1 \pmod{3}$ be a prime such that $2$ is a cube
modulo $p$. Then
\[
\dim_{\F_2} \Pic^0(C_p)/2\Pic^0(C_p)=2.
\]
Moreover, each of the $4$ elements of $\Pic^0(C_p)/2 \Pic^0(C_p)$ can be represented by
the class of $P-\OO$ for some 
$P \in C(\F_p)$. 
\end{lem}
\begin{proof}
We can put $E=(C,\OO)$ in Weierstrass form by sending
$\OO$ to the point at infinity. A Weierstrass model is
\[
y^2 + y = x^3 - 7. 
\]
The $2$-division polynomial of this model is $4x^3-27$. 
The assumptions on $p$ ensure that the $2$-division polynomial
splits completely over $\F_p$; thus $E$ has
full $2$-torsion over $\F_p$. 
\end{proof}

\section{Reduction of Lines}\label{sec:red}
In what follows we would like to conveniently parametrize a $\Q_p$-line $\ell$
in $\PP^3$. Of course if $P$ and $Q$ are two distinct points on $\ell$
then $\ell(\Q_p)=\{ s P+ t Q : (s:t) \in \PP^1(\Q_p)\}$. Now if $S$
is a cubic surface given by $F=0$ then the points of $\ell \cdot S$
correspond to the roots of $F(sP+tQ)$ with the correct
multiplicity. The problem with such a parametrization is that
it is possible that $\overline{P}=\overline{Q}$ and so we do not
obtain a parametrization of $\overline{\ell}$ by reducing the
parametrization of $\ell$. In this brief section we indicate how to make
a good choice of parametrization. 
It will be convenient to use the identification of lines in $\PP^3$
with planes in $4$-dimensional space passing through the origin.
Let $V_\ell$ be the
$2$-dimensional $\Q_p$-subspace of $\Q_p^4$ generated by the points
of $\ell(\Q_p)$. Let $W_\ell=V_\ell \cap \Z_p^4$.
This $W_\ell$
is a $\Z_p$-module of rank $2$ and we let $\uu$, $\vv$ be a $\Z_p$-basis.
Then the line $\ell$ can be parametrized as $s \uu+t \vv$, and
$\overline{\ell}$ as $s \overline{\uu}+t \overline{\vv}$.
We shall call this a {\em good} parametrization for $\ell$.
We note that if $P\in \ell(\Q_p)$ then there is a coprime
pair of $p$-adic integers $\lambda$, $\mu$ such that
$P=\lambda \uu+\mu\vv$; this can be easily seen by regarding
$P$ as a primitive element of $W_\ell$.

We now turn our attention to the surfaces of Theorems~\ref{thm:two}
and~\ref{thm:three}. Let $M$ denote a non-zero squarefree integer
and let $S$ be either of the surfaces in \eqref{eqn:one}
or \eqref{eqn:two}. 
Let $p \ne 3$ be a prime dividing $M$;
we would like to study the group $H_S(\Q_p)$.
Denote the reduction of $S$ modulo $p$ 
by $\overline{S}$. It is seen that $\overline{S}$ is a cone
over $C_p=C \times \F_p$, where $C$ is given in \eqref{eqn:C},
with vertex at $(\overline{0}: \overline{0} : \overline{0} : \overline{1})$.
It is convenient to split the set of $p$-adic points on $S$ 
into subsets of bad and good reduction:
\[
\Sb=\{P \in S(\Q_p) : \overline{P}=(\overline{0}:\overline{0}:\overline{0}:\overline{1})\},
\qquad
\Sg=\{P \in S(\Q_p) : \overline{P} \neq (\overline{0}:\overline{0}:\overline{0}:\overline{1})\}.
\]
We note however that $\Sb=\emptyset$ 
in the case of the surface $S_M^\prime$ of \eqref{eqn:two}. Now
we can think of points on $\Sg$ as reducing to points on $C(\F_p)$.
We define $\phi : \Sg \rightarrow C(\F_p)$ by 
$\phi(x:y:z:w)=(\overline{x}:\overline{y}:\overline{z})$ where the four
coordinates $x$, $y$, $z$, $w$ are taken to be coprime in $\Z_p$.

\begin{lem}\label{lem:Sdotlred}
Let $M$ denote a non-zero squarefree integer and $S$
either of the surfaces in \eqref{eqn:one}
or \eqref{eqn:two}. Let $p \neq 3$ be a prime dividing $M$.
Let $\ell$ be a $\Q_p$-line such that $\ell \cdot S=P_1+P_2+P_3$
with $P_i \in \Sg$. Suppose $\overline{\ell} \not \subset \overline{S}$.
Then $\phi(P_1)+\phi(P_2)+\phi(P_3) \sim 3 \overline{\OO}$ in $\Pic(C_p)$.
\end{lem}
\begin{proof}
Let $s \uu+t\vv$ be a good parametrization of $\ell$ in
the above sense.
There are coprime pairs $\lambda_i$, $\mu_i \in \Z_p$ for $i=1$, $2$, $3$
such that
$P_i=\lambda_i \uu + \mu_i \vv$ for $i=1$, $2$, $3$.
Since $\ell \cdot S=P_1+P_2+P_3$ we have
\begin{equation}\label{eqn:Sdotl}
F(s\uu+t\vv)=\alpha \prod_{i=1}^3(\mu_i s -\lambda_i t)
\end{equation}
for some $\alpha$ in $\Z_p$. 
Moreover $\alpha \not \equiv 0 \pmod{p}$ as $\overline{\ell} \not \subset \overline{S}$. 
Write $\uu=(u_0,u_1,u_2,u_3)$ and $\vv=(v_0,v_1,v_2,v_3)$
and let 
$\uu^\prime=(\overline{u_0},\overline{u_1},\overline{u_2})$ 
and $\vv^\prime=(\overline{v_0},\overline{v_1},\overline{v_2})$. 
We first show that $\uu^\prime$ and $\vv^\prime$ are 
$\F_p$-linearly independent. If not, then $\overline{\ell}$ passes through
the vertex $(\overline{0}:\overline{0}:\overline{0}:\overline{1})$
as well as the $\overline{P_i}$ which are distinct from the vertex.
This forces $\overline{\ell} \subset \overline{S}$ contradicting
the lemma's assumption that $\overline{\ell} \not \subset \overline{S}$. Thus $\uu^\prime$
and $\vv^\prime$ are $\F_p$-linearly independent.

Let $\ell^\prime \subset \PP^2$ 
be the $\F_p$ line $s \uu^\prime +t \vv^\prime$. 
Reducing \eqref{eqn:Sdotl} modulo $p$ we instantly see that
$\ell^\prime \cdot C_p=\phi(P_1)+\phi(P_2)+\phi(P_3)$
which concludes the proof of the lemma.
\end{proof}

\section{A Homomorphism for $H_{S_M}(\Q_p)$}\label{sec:SM}
In this section $M$ denotes a non-zero squarefree integer,
and $S$ the cubic surface denoted by $S_M$ in \eqref{eqn:one}:
\[
S : x^3+y^3+z(z^2+M w^2)=0.
\]
Let $C$ be the plane cubic curve \eqref{eqn:C}.
Again, let $p$ be a prime divisor of $M$ different from $3$ and let $C_p=C \times \F_p$. 
In this section we shall define a 
surjective homomorphism $H_S(\Q_p) \rightarrow \Pic^0(C_p)/2\Pic^0(C_p)$.
As previously observed, the reduction of $S$ modulo $p$ is a cone over $C_p$ with vertex at 
$ (\overline{0}:\overline{0}:\overline{0}:\overline{1})$.
What makes the situation here somewhat tricky is that  
this singular point lifts to $p$-adic (and even rational) points on $S$; for example $(0:0:0:1) \in S(\Q)$
is such a lift. Thus $\Sg$ is strictly smaller than $S(\Q_p)$.
We shall extend $\phi:\Sg \rightarrow C(\F_p)$ defined in Section~\ref{sec:red}
to $\phi:S(\Q_p) \rightarrow C(\F_p)$ where $\phi(P)=\overline{\OO}$
for any $P \in \Sb$.

\begin{prop}\label{prop:SM}
Let $\psi : S(\Q_p) \rightarrow \Pic^0(C_p)$ be given by $\psi(P)=\phi(P)-\overline{\OO}$.
Then $\psi$ induces a well-defined surjective homomorphism
\[
\psi : H_S(\Q_p) \rightarrow \Pic^0(C_p)/2 \Pic^0(C_p),
\qquad \psi([P])=\psi(P) \pmod{2 \Pic^0(C_p)}.
\]
\end{prop}

Before proving Proposition~\ref{prop:SM} we shall need the following 
three lemmas.

\begin{lem}\label{lem:zzero}
Suppose $P\in S(\Q_p)$ such that $\overline{P}$
lies in the plane $z=\overline{0}$. Then 
$\psi(P) \in 2 \Pic^0(C_p)$.
\end{lem}
\begin{proof}
If $P \in \Sb$ then $\psi(P)=0$. Otherwise, $\phi(P)$
is a flex point on $C_p$. Hence $\psi(P)=\phi(P)-\overline{\OO}$
is a difference of two flexes and so is an element of order 
dividing $3$ in $\Pic^0(C_p)$.
It follows that $\psi(P) \in 2 \Pic^0(C_p)$:
indeed $\psi(P)=-2\psi(P)$.
\end{proof}

\begin{lem}\label{lem:badline}
Let $\ell$ be a $\Q_p$-line not contained in $S$. Suppose that
$\ell \cdot S=P_1+P_2+P_3$ where $P_i \in S(\Q_p)$.
Suppose moreover that $\overline{\ell} \subset \overline{S}$.
Then $\psi_1(P)+\psi_2(P)+\psi_3(P) \in 2\Pic^0(C_p)$.
\end{lem}
\begin{proof}
Let $s \uu+t\vv$ be a good parametrization of $\ell$.
In particular we know that $\uu$, $\vv \in \Z_p^4$
and $\overline{\uu}$, $\overline{\vv}$ are independent
modulo $p$, and that there are coprime pairs $\lambda_i$,
$\mu_i$ such that $P_i=\lambda_i \uu+\mu_i \vv$.
Now $\overline{\ell}$ is contained in $\overline{S}$
and so must pass through the
vertex $(\overline{0}:\overline{0} : \overline{0} : \overline{1})$.
Thus, applying a unimodular transformation, we may assume
that our pair $\uu$, $\vv$ have the form
\begin{equation}\label{eqn:uv}
\uu=(u_1,u_2,u_3,0), \qquad \vv=(p v_1, p v_2, p v_3, 1), \qquad
u_i,~v_i \in \Z_p.
\end{equation}
Let $F=x^3+y^3+z^3+Mz w^2$. We can also assume that $F(\vv) \ne 0$,
by replacing $\vv$ with $\vv+\alpha p \uu$ for an appropriate $\alpha \in \Z_p$.
It follows that the polynomial $F(\uu+t \vv)$ is of degree $3$
and has
precisely $3$ roots in $\Q_p$: namely $\mu_i/\lambda_i$ for $i=1,2,3$.
However, expanding $F(\uu+t\vv)$ we obtain
\begin{equation}\label{eqn:F}
\sum u_i^3+3p\left(\sum u_i^2 v_i\right) t+
\left(3 p^2 \sum u_i v_i^2+M u_3\right)t^2+
\left(p^3 \sum v_i^3+pM v_3\right)t^3. 
\end{equation}
Write $F(\uu+t\vv)=a_0+a_1 t +a_2 t^2+a_3 t^3$, and let
 $\alpha_i=\ord_p(a_i)$.
Since $\overline{\ell} \subset \overline{S}$ we see that $\alpha_i \geq 1$.
Suppose first that $p \nmid u_3$. 
As $M$ is squarefree and $p \mid M$, we see that
\[
\alpha_0 \geq 1, \qquad \alpha_1 \geq 1,
\qquad \alpha_2=1, \qquad \alpha_3 \geq 2.
\]
We shall need to study the Newton polygon of this polynomial
which is the convex hull of the four points $(i,\alpha_i)$;
see for example \cite[page 19]{Koblitz}.
The Newton polygon  contains precisely 
one segment of positive slope joining
$(2,1)$ with $(3,\alpha_3)$. The other
segments have non-positive slope. 
Hence the polynomial has precisely
one root with negative valuation and
two roots in $\Z_p$. By reordering the
$P_i$ we may suppose that $\mu_1/\lambda_1$
has negative valuation and $\mu_i/\lambda_i \in \Z_p$
for $i=2$, $3$. From the expressions $P_i=\lambda_i \uu+\mu_i \vv$
we see that $P_1 \in \Sb$ and $P_2$, $P_3 \in \Sg$.
Moreover $\phi(P_2)=\phi(P_3)=(u_1 : u_2 :u_3) \in C(\F_p)$.
Hence $\psi(P_1)+\psi(P_2)+\psi(P_3) \in 2 \Pic^0(C_p)$.

Finally we must deal with the case
$p \mid u_3$. In this case $\overline{P}_i$
lie in the plane $z=\overline{0}$.
The lemma follows from Lemma~\ref{lem:zzero}.
\end{proof}

\begin{lem}\label{lem:ol}
Let $\ell$ be a $\Q_p$-line contained in $S$. For every $P\in \ell(\Q_p)$,
we have $\psi(P) \in 2 \Pic^0(C_p)$.
\end{lem}
\begin{proof}
Since $\overline{\ell}$
passes through the vertex of $\overline{S}$ we may parametrize
as $s \uu+t \vv$ where $\uu$, $\vv$ are as in \eqref{eqn:uv}.
Now the polynomial in \eqref{eqn:F} vanishes identically.
From the coefficient of $t^2$ we see that $p \mid u_3$.
Hence $\overline{\ell}$ lies in the plane $z=\overline{0}$.
Now the lemma follows from Lemma~\ref{lem:zzero}.

In fact, more is true:
the only $\Q_p$-lines on $S$ are contained in the $z=0$
plane, but we do not need this.
\end{proof}
\begin{proof}[Proof of Proposition~\ref{prop:SM}]
Hensel's Lemma shows that the map $\phi: S^{\mathrm{gd}} \rightarrow C(\F_p)$
is surjective, and thus $\psi$ is surjective.
Therefore,
it is enough to show that $\psi(P_1)+\psi(P_2)+\psi(P_3) \in 2 \Pic^0(C_p)$
whenever $P_1$, $P_2$, $P_3 \in S(\Q_p)$ and
\begin{enumerate}
\item[(i)] either there is a $\Q_p$-line $\ell$ not contained
in $S$ with $\ell \cdot S=P_1+P_2+P_3$,
\item[(ii)] or there is a $\Q_p$-line $\ell$
contained in $S$ with $P_1$, $P_2$, $P_3 \in \ell$. 
\end{enumerate}
For (ii) we know by Lemma~\ref{lem:ol} that $\psi(P_i) \in 2 \Pic^0(C_p)$,
so it remains to deal with (i).
Thus suppose that $\ell$ is a $\Q_p$-line not contained in $S$
and $\ell \cdot S=P_1+P_2+P_3$ where $P_i \in S(\Q_p)$.
If all three $P_i \in \Sb$ then the required
result follows from the definition of $\psi$.
Suppose that at least one $P_i\in \Sg$. 
\begin{itemize}
\item 
If all three $P_i \in \Sg$ and $\overline{\ell}
\not \subset \overline{S}$ then we can conclude using Lemma~\ref{lem:Sdotlred}.
\item If all three $P_i \in \Sg$  and $\overline{\ell} \subset \overline{S}$ 
we can conclude using Lemma~\ref{lem:badline}.
\end{itemize}
We have reduced to the case where at least one of the $P_i$ is in $\Sg$
and at least one is in $\Sb$. This forces $\overline{\ell} \subset \overline{S}$
and again we can conclude using Lemma~\ref{lem:badline}.
\end{proof}

\section{Proof of Theorem~\ref{thm:two}}\label{sec:proofone}
In this section we shall put together the results of the
previous sections 
to prove Theorem~\ref{thm:two}. Thus
let $p_1,\dots,p_s$ ($s \geq 1$) be primes such that
$p_i \equiv 1 \pmod{3}$ and $2$ is a cube modulo each $p_i$.
Let $M=\prod p_i$ and $S=S_M$ be the cubic surface in \eqref{eqn:one}.

For $p=p_i$,  Proposition~\ref{prop:SM} gives a surjective homomorphism
\[
\psi : H_S(\Q_{p}) \rightarrow \Pic^0(C_{p})/2\Pic^0(C_{p}).
\]
Now $\psi(1:-1:0:0)=0$, thus the restriction of $\psi$
to $H_S^0(\Q_p)$ is still surjective. In particular,
\[
\dim_{\F_2} H^0_S(\Q_p)/2 H^0_S(\Q_p) \ge \dim_{\F_2} \Pic^0(C_{p})/2\Pic^0(C_{p})=2
\]
by Lemma~\ref{lem:mod2}. However, by Proposition~\ref{prop:dense}, the map
\[
H^0_S(\Q)/2H^0_S(\Q) \rightarrow \prod_{p=p_1}^{p=p_s} H^0_S(\Q_p)/2 H^0_S(\Q_p)
\]
is surjective, and so
\[ 
\dim_{\F_2} H^0_S(\Q)/2 H^0_S(\Q) \ge 2s.
\]
As $S$ contains the $\Q$-rational line
$x+y=z=0$, we know by Theorem~\ref{thm:line} that
$H^0_S(\Q) = H_S(\Q)[2]$. Thus
$\dim_{\F_2} H_S(\Q)[2] \geq 2s$. Now the point
$(1:-1:0:0)$ is $\Q$-ternary on $S$.
Thus we can apply Theorem~\ref{thm:terngen} obtaining
$r(S,\Q)\geq \dim_{\F_2} H_S(\Q)[2] \geq 2s$, which proves
Theorem~\ref{thm:two}.

\section{A Homomorphism for $H_{S^\prime_M}(\Q_p)$}\label{sec:SMd}

In this section $M>1$ denotes a squarefree integer
and $S$ the cubic surface denoted by $S^\prime_M$ in \eqref{eqn:two}:
\[
S : x^3+y^3+z^3+Mw^3=0.
\]
Let $C$ be the plane cubic curve \eqref{eqn:C}.
Let $p$ be a prime divisor of $M$ different from $3$
and 
  $\phi : S(\Q_p) \rightarrow C(\F_p)$
as defined in Section~\ref{sec:red}; as previously observed for this surface,
$\Sb=\emptyset$ and so $\Sg=S(\Q_p)$. 
Let $\psi : S(\Q_p) \rightarrow \Pic^0(C_p)$ be given by
$\psi(P)=\phi(P)-\overline{\OO}$.

\begin{prop}\label{prop:SMd}
Let $C_p=C \times \F_p$. Then $\psi$ extends uniquely to
a well-defined surjective homomorphism
\[
\psi : H_S(\Q_p) \rightarrow \Pic^0(C_p)/3 \Pic^0(C_p),
\qquad \psi([P])=\psi(P) \pmod{3 \Pic(C_p)}.
\]
\end{prop}
\begin{proof}
Surjectivity follows as in the proof of Proposition~\ref{prop:SM}.
It is therefore 
enough to show that $\psi(P_1)+\psi(P_2)+\psi(P_3) \in 3 \Pic^0(C_p)$
whenever $P_1$, $P_2$, $P_3 \in S(\Q_p)$ and
\begin{enumerate}
\item[(i)] either there is a $\Q_p$-line $\ell$ not contained
in $S$ with $\ell \cdot S=P_1+P_2+P_3$,
\item[(ii)] or there is a $\Q_p$-line $\ell$
contained in $S$ with $P_1$, $P_2$, $P_3 \in \ell$. 
\end{enumerate}
The second possibility does not arise since none of the $27$ lines
on $S$ are defined over $\Q_p$. Thus suppose $\ell$ is a
$\Q_p$-line not contained in $S$ such that
$\ell \cdot S=P_1+P_2+P_3$. If $\overline{\ell} \not \subset \overline{S}$
then by Lemma~\ref{lem:Sdotlred} we know that
$\phi(P_1)+\phi(P_2)+\phi(P_3) \sim 3 \overline{\OO}$ in $\Pic(C_p)$.
In this case $\psi(P_1)+\psi(P_2)+\psi(P_3)=0$ in $\Pic^0(C_p)$.
Thus we may suppose that $\overline{\ell} \subset \overline{S}$.
Hence $\overline{\ell}$ passes through the 
vertex $(\overline{0}:\overline{0}:\overline{0}:1)$, 
otherwise its projection on to the $xyz$-projective plane
would be a line contained in irreducible curve $C_p$ which is impossible. 
It follows that
$\phi(P_1)=\phi(P_2)=\phi(P_3)$ and thus $\psi(P_1)+\psi(P_2)+\psi(P_3) \in 3 \Pic^0(C_p)$. 
\end{proof}

\section{Proof of Theorem~\ref{thm:three}}\label{sec:prooftwo}
This follows from Lemma~\ref{lem:mod3} and
  Propositions~\ref{prop:densed} and~\ref{prop:SMd}
in exactly the same way as the proof of Theorem~\ref{thm:two} 
(Section~\ref{sec:proofone}).

\end{document}